\documentclass[a4paper,12pt,onecolumn]{article}

\usepackage{hyperref}
\usepackage[vmargin=2cm,hmargin=2cm,headheight=14.5pt,top=2cm,headsep=.5cm]{geometry}
\usepackage[utf8]{inputenc}
\usepackage{bm}
\usepackage{empheq}

\usepackage{stackrel}
\usepackage{cases}
\usepackage{mathtools}
\usepackage{amsthm,amsmath,amscd}
\usepackage{makeidx}
\usepackage[charter]{mathdesign}
\usepackage{tikz-cd}
\tikzcdset{every label/.append style = {font = \small}}

\usepackage{graphicx}
\DeclareMathSizes{12}{12}{8}{6}

\usepackage{cite}
\usepackage{url}
\usepackage[ddmmyy]{datetime}

\usepackage{marvosym}
\usepackage{perpage}
\MakePerPage[1]{footnote}

\newtheoremstyle{ptheorem}{1em}{0em}{\itshape}{}{\bfseries}{.}{.5em}{\thmname{#1}\thmnumber{ #2}\thmnote{ (\hspace{-.01pt}{#3})}}

\theoremstyle{ptheorem}

\newtheorem{thm}{Theorem}[section]

\newtheorem{con}[thm]{Conjecture}

\newtheoremstyle{hdef}{1em}{0em}{}{}{\bfseries}{.}{.5em}{\thmname{#1}\thmnumber{ #2}\thmnote{ (\hspace{-.01pt}{#3})}}
\theoremstyle{hdef}

\makeatletter
\newtheoremstyle{premark}{1em}{0em}{
\addtolength{\@totalleftmargin}{1.5em}
\addtolength{\linewidth}{-1.5em}
\parshape 1 1.5em \linewidth}{}{\scshape}{.}{.5em}{}
\makeatother

\theoremstyle{premark}

\numberwithin{equation}{section}
\numberwithin{figure}{section}

\DeclareMathOperator{\tr}{tr}

\newcommand{\cC}{{\mathcal C}}

\newcommand{\cM}{{\mathcal M}}

\newcommand{\bC}{{\mathbb C}}

\newcommand{\bN}{{\mathbb N}}

\newcommand{\bR}{{\mathbb R}}

\renewcommand{\phi}{\varphi}

\newcommand{\ol}{\overline}

\title{Differential systems with reflection\\ and matrix invariants}

\author{Santiago Codesido\footnotemark \\ F. Adri\'an F. Tojo\footnotemark}

\date{\today}
\begin{document}
 \maketitle
 
%

%
\footnotetext{\footnotemark Département de Physique Théorique et Section de Mathématiques. Université de Genève, Genève, CH-1211 Switzerland. \href{mailto:santiago.codesido@unige.ch}{santiago.codesido@unige.ch}}
\footnotetext{Corresponding author. Instituto de Matemáticas, Universidade de Santiago de Compostela,  15782, Facultade de Matemáticas, Santiago, Spain. \href{mailto:fernandoadrian.fernandez@usc.es}{fernandoadrian.fernandez@usc.es}. Partially supported by project MTM2016-75140-P (AEI/FEDER, UE) and Xunta de Galicia (Spain), project EM2014/032.}

\medbreak

\begin{abstract}
In this work we derive important properties regarding matrix invariants which occur in the theory of differential equations with reflection.
\end{abstract}

\textbf{Keywords:} differential equations with reflection, matrix invariants.

\section{Introduction}

In recent works regarding the solution and Green's functions of Differential Equations with Reflection (see for instance \cite{Toj3, TojMI, CTMal,CT17}) the strong relation between linear analysis and linear algebra is highlighted. In particular, in the most recent of the aforementioned works, the authors obtain an explicit fundamental matrix for the system of differential equations with reflection
\begin{equation}\label{hlsystem}Hu(t):=Fu'(t)+Gu'(-t)+A u(t)+Bu(-t)=0, t\in\bR,
\end{equation}
where $n\in\bN$, $A,B,F,G\in\cM_n(\bR)$ and $u:\bR\to\bR^n$. To be precise, they prove the following result.
\begin{thm}[\cite{CT17}]\label{thmexpfm}
	Assume $F-G$ and $F+G$ are invertible. Then
	\begin{equation*}\label{Xseries}X(t): =     \sum_{k=0}^\infty\frac{E^k t^{2k}}{(2k)!}  -(F+G)^{-1}(A+B)\sum_{k=0}^\infty\frac{E^k t^{2k+1}}{(2k+1)!},\end{equation*}
	where $E=(F-G)^{-1}(A-B)(F+G)^{-1}(A+B)$, is a fundamental matrix of problem \eqref{hlsystem}. If we further assume $A-B$ and $A+B$ are invertible, then $E$ is invertible and we can consider a square root $\Omega$  of $E$. Then,
	\begin{equation*}\label{fme}X(t)=\cosh \Omega t -(F+G)^{-1}(A+B)\Omega^{-1}\sinh\Omega t.\end{equation*}
\end{thm}

What is more, in another recent work  the authors proved an analog of the Liouville's formula for the case with reflections in systems of order two.

\begin{thm}[Abel-Jacobi-Liouville Identity \cite{CoTo17}] \label{AJLid}Let $n=2$ in equation \eqref{hlsystem}. Then $(|X|,|X'|)$ is the unique solution of the system of differential equations
	\begin{equation*}\label{n2e}\begin{aligned}x''= &\tr(E)x-2y, \\y''= &-2|E|x+\tr(E)y,\end{aligned}
	\end{equation*}
	subject to the one point conditions \[x(0)=1,\quad y(0)=|M_+|,\quad x'(0)=-\tr(M_+),\quad y'(0)=\tr(\operatorname{Adj}(M_+)E).\]
\end{thm}

The authors also presented in that work the following conjecture: 
\begin{con} For any $n\ge 1$, if $X(t)$ is a fundamental matrix of problem \eqref{hlsystem}, then $|X(t)|$ can be obtained as a component of the solution of a linear system of differential equations with constant coefficients, those coefficients depending only on the different matrix invariants of $E$, which is defined as in Theorem \ref{thmexpfm}.
\end{con}

In order to attempt proving this conjecture, and taking into account the proof of Theorem \ref{AJLid}, we need to study the different matrix invariants of the matrices appearing in the theory.
\section{The $\mathbf Y$ matrix}
For $X(t)$ the fundamental matrix of the problem, define 
\begin{equation}
Y(t) := X(t)^{-1} X'(t).
\label{Ydef}
\end{equation}
We have that $X = S_1-M_+ S_2$ where $S_1$ and $S_2$ s are power series in $E$ which we can formally give, by using $\Omega = \sqrt{E}$, as $S_1=\cosh(\Omega  t),\    S_2=\Omega^{-1} \sinh(\Omega t).$

Notice both are indeed power series in $\Omega^2 = E$. Since $X'' = X E$ and $   (Y^{-1})'= - Y^{-1} Y'Y^{-1}$ we have that
\begin{equation}
Y'= E - Y^2.
\label{ricattiEq}
\end{equation}
Using the construction of \cite{levin} we build an associated ODE system
\begin{equation*}
z' = \left(\begin{array}{cc}
0 & E \\ I & 0
\end{array}\right)z.
\end{equation*}
The system has as fundamental matrix
\begin{equation*}
\left(\begin{array}{cc}
\cosh(\Omega t) & \Omega^{-1}\sinh(\Omega t) \\ \Omega \sinh(\Omega t) & \cosh(\Omega t)
\end{array}\right).
\end{equation*}

The solution of equation (\ref{ricattiEq}) is then given by
\begin{equation*}
Y(t) = \left[\cosh(\Omega   t) Y(0)+\Omega \sinh(\Omega t) \right] \left[\Omega^{-1}\sinh(\Omega   t) Y(0)+\cosh(\Omega t) \right]^{-1},
\end{equation*}
which in terms of the $S$ functions is
$Y(t) = \left[-S_1 M_+ +E S_2\right] \left[-S_2 M_+ + S_1\right]^{-1}$,
where we fix the initial condition with $Y(0)=X'(0)=-M_+$.
This seems like a commuted version of expression (\ref{Ydef}), but it is nothing more than the hypergeometric identity.


Consider the Liouville equation for $Y$ itself, that is,
\[
\left(\log |Y|\right) '= \mathrm{Tr}\left(Y^{-1} Y'\right) = \mathrm{Tr}\left( Y - Y^{-1}E \right).\] Then we have
$\mathrm{Tr}(Y^{-1}E)=\mathrm{Tr}(Y)-\left(\log |Y|\right)'$, which can be calculated in terms of invariants of $|X|$.
	
\section{Complex systems}
The main involution occurring in the theory of complex variable is the complex conjugation $\cC:\bC\to \bC$, $\cC(z)=\ol z$. It is, in fact, a reflection with respect to the second variable if we write $z=(x,y)\in\bR^2$: $\cC(x,y)=(x,-y)$.

We consider now an operator $L$ acting on $z(t)$ as
\begin{equation}
\label{mainComplexEquation}
A_0 z(t) + A_1 \overline{z(t)} + B_0 z(t)+ B_1 \overline{z'(t)},
\end{equation}
where $z: \mathbb{R}\to \mathbb{C}^n$, and $A_i,B_i \in \mathcal{M}_{n\times n}\left(\mathbb{C}\right):=\mathcal{M}$.

We can consider a extended algebra $\mathcal{M}^*$ by with the linear operation of complex conjugation which acts as
$
\mathcal{C} z = \overline{z},   z\in \mathbb{C}^n.
$

It is easy to see that the following properties hold\footnote{In fact, one could use here any involution for which matrix conjugation verifies $\mathcal{C}A\mathcal{C}\in \mathcal{M}$
	by defining $\overline{A}$ suitably.},
\begin{equation}
\label{conjugationProperties}
\mathcal{C}^2 = I,    \mathcal{C} A = \overline{A} \mathcal{C},
\end{equation}
where $\overline{A}$ is the complex conjugate of $A$.

The we consider the free product quotiented by these relations,
\begin{equation*}
\mathcal{M}^* = \mathcal{M} \star \left\{\mathcal{C}\right\} \left/ \left( \mathcal{C}^2 = I, \mathcal{C} A = \overline{A} \mathcal{C}\right).\right.
\end{equation*}
Now, we can see that this is in fact a $\mathbb{Z}_2$-graded algebra. Due to the conditions (\ref{conjugationProperties}), we can move any $\mathcal{C}$'s to the
right, and any power of it is reduced modulo $2$. Therefore, any element of $A\in \mathcal{M}^*$ can be written as
$
\mathbf{A} = A_0 + A_1 \mathcal{C}$ with $A_0,A_1 \in \mathcal{M}
$.
As a vector space, $
\mathcal{M}^* = \mathcal{M} \oplus \mathcal{M}$.

The grading is clear by looking at the product of two generic elements,
\begin{equation}
\label{conjugationAlgebraProduct}
\mathbf{A} \mathbf{B}=\left(A_0 + A_1\mathcal{C}\right)\left(B_0 + B_1\mathcal{C}\right) =\left(A_0 B_0 + A_1 \overline{B_1}\right) + \left( A_0 B_1 + A_1 \overline{B_0} \right) \mathcal{C}.
\end{equation}
We can also calculate a explicit inverse, using $
\left(I-A \mathcal{C}\right)\left(I+A \mathcal{C}\right) = I - A \overline{A}
$, 
from where
$\left(I+A \mathcal{C}\right)^{-1} = \left(I-A \overline{A}\right) \left(I - A \mathcal{C}\right)$.

Since $\left(A B\right)^{-1}=B^{-1} A^{-1}$, we can generalize it to
\begin{equation*}
\left(A_0 + A_1 \mathcal{C}\right)^{-1} = \left(A_1^{-1}A_0 - \overline{A_0^{-1}A_1}\right)^{-1} \left(A_1^{-1}-\overline{A_0^{-1}}\mathcal{C}\right).
\end{equation*}
As it is, the expression is unclear when either $A_i=0$. We rewrite
\begin{equation}
\label{conjugationAlgebraInverse}
\mathbf{A}^{-1} = \Delta\left(A_0,A_1\right) + \Delta\left(\overline{A_1},\overline{A_0}\right) \mathcal{C},
\end{equation}
with
\begin{equation*}
\Delta\left(A_0,A_1\right)=\begin{cases}
0, & A_0 = 0, A_1 \neq 0, \\
\left(A_0 - A_1 \overline{A_0^{-1}A_1}\right)^{-1}, & A_0 \neq 0, A_1 \neq 0.
\end{cases}
\end{equation*}

As a last note, this can of course be realized as a matrix algebra over $\mathbb{R}^{2n}$, although it does not lead to anything new (other than clutter). For the record,
one can take a representation
\[
\rho\left(z\right) = \left(\begin{array}{c}
\Re z \\ \hline \Im z
\end{array}\right), \quad
\rho\left(A\right) = \left(\begin{array}{c|c}
\Re A & -\Im A \\ \hline
\Im A & \Re A
\end{array}\right), \quad
\rho\left(\mathcal{C}\right) = \left(\begin{array}{c|c}
I & 0 \\ \hline 0 & -I
\end{array}\right),
\]
for which it is easy to see that properties such as (\ref{conjugationProperties}) or (\ref{conjugationAlgebraProduct}) hold.

\subsection{Equation reduction}	
With these tools, one can rewrite equation (\ref{mainComplexEquation}) as
$
\mathbf{B} z'(t) + \mathbf{A} z(t)  = 0,
$
which can be reduced to
$
z'(t) + \left(\mathbf{B}^{-1} \mathbf{A}\right) z(t) =0
$. 
This means we can just focus on the study of
\begin{equation*}
z'(t) + A_0 z(t) + A_1 \overline{z(t)} =0.
\end{equation*}

Of course, one could now look for a fundamental operator inside the $\mathcal{M}^*$ algebra, such that solutions fulfill
\begin{equation*}
z(t) = \mathbf{X}(t) z(0),
\end{equation*}
which is
\begin{equation*}
\mathbf{X} = \cosh \left(\mathbf{A} t\right) - \sinh \left(\mathbf{A} t\right).
\end{equation*}
Unfortunately, it does not seem easy to calculate terms like $\mathbf{A}^n$ at the moment. We could in principle directly get
an explicit fundamental matrix if we could write a manageable expression. However, we do learn something important. Since $\mathbf{A}^n \in \mathcal{M}^*$, then
$\cosh\left( \mathbf{A} t\right) \in \mathcal{M}^*$ too. Hence, if we want to find fundamental matrices for the problem, the ansatz must be of the form
\begin{equation*}
z(t) = \left( X_0(t) + X_1(t) \mathcal{C} \right)z(0)= X_0(t) z(0) + X_1(t) \overline{z(0)},
\end{equation*}
where $X(0)=I$ and $Y(0)=0$ as to agree with $\mathbf{X}(0)=I$.

\subsection{$\mathbf{A^n}$ generating function}	
The components of $\mathbf{A}^n$ can still be algorithmically computed  by using the expression
\begin{equation*}
\left(I-t \mathbf{A}\right)^{-1} = \sum_{n=0}^\infty \left(t \mathbf{A}\right).
\end{equation*}
By writing the explicit inverse of $I-t\mathbf{A}=\left(I-t A_0\right)-t A_1 \mathcal{C}$, we get
\begin{small}
	\begin{equation*}
	\mathbf{A}^n = \frac{1}{n!} \frac{\mathrm{d}}{\mathrm{d}t}\left[ \left(A_1^{-1}\left( I - t A_0\right) - \overline{\left(I -t A_0\right)^{-1}A_1}\right)^{-1}
	\left(A_1^{-1} - \overline{ \left(t^{-1}I-A_0\right)^{-1} } \mathcal{C}\right)\right]_{t=0}.
	\end{equation*}
\end{small}

\subsection{Ansatz}
We take the system
\begin{equation*}
z'+A z+\overline{B z} = 0,    z: \mathbb{R}\to \mathbb{C}^n,   A,B \in \mathcal{M}
\end{equation*}
and introduce the ansatz
\begin{equation*}
z = X z_0 + \overline{Y z_0},    z_0=z(0),    X,Y: \mathbb{R} \to \mathcal{M},
\end{equation*}
\begin{equation*}
\left(X' + A X + \overline{B} Y \right) z_0 + \left(Y'+ \overline{A} Y + B X\right) \overline{z_0} = 0,
\end{equation*}
\begin{equation}
\label{fundamentalSystem}
\left\{ \begin{array}{r}
X' + A X + \overline{B} Y = 0, \\
Y'+ \overline{A} Y + B X = 0.
\end{array} \right.
\end{equation}
This a ordinary system. Take $X''$ and substitute $Y'$ and $Y$ through equation (\ref{fundamentalSystem}),
\begin{align*}
\begin{split}
& X'' + A X' + \overline{B} \left(-B X - \overline{A}   \overline{B^{-1}}\left(-X'-A X\right)\right) = 0, \\
& X'' + \left(A + \overline{B A B^{-1}}\right) X' + \left(\overline{B A B^{-1}}A - \overline{B}B\right) X = 0.
\end{split}
\end{align*}
Unsurprisingly, we get a very similar structure to the inverses in expression (\ref{conjugationAlgebraInverse}).For the sake of notation, we will rename the coefficients as
\begin{equation*}
\label{singleFundamentalMatrixEquation}
X'' + F X' + G X= 0.
\end{equation*}
Repeating the process, we get 
\begin{equation*}
Y'' + \overline{F} Y' + \overline{G} Y = 0.
\end{equation*}
The initial conditions for this second order problem are given by $\mathbf{X}(0)=I$ and equation (\ref{fundamentalSystem}). That is,
\[X(0) = I, \quad Y(0)=0,\quad
X'(0) = -A, \quad Y'(0) = -\overline{A}.
\]

In principle, we could now take as an ansatz
\begin{equation*}
X = \alpha e^{(\Gamma + \Omega) t} + \beta e^{(\Gamma-\Omega) t},
\end{equation*}
subject to the conditions
\begin{equation*}
X(0) = \alpha + \beta,    X'(0) = \alpha (\Gamma+\Omega) + \beta (\Gamma-\Omega),
\end{equation*}
which can be inverted into
\begin{align*}
\begin{split}
\alpha = &\frac{1}{2} \left[ X'(0) - X(0) \left(\Gamma-\Omega\right) \right] \Omega^{-1},\\
\beta = &-\frac{1}{2} \left[ X'(0) - X(0) \left(\Gamma+\Omega\right) \right] \Omega^{-1} .
\end{split}
\end{align*}

\section{Generalized matrix invariants}
%
%

In the following section we use the concept of crossed or generalized matrix invariants which can be found in \cite{simon} and \cite{cgm} among others.
\subsection{Definition and basic properties}
Let $X_1,\dots,X_N \in GL\left(n\right)$. Define

\begin{equation}\label{z}	Z\left(X_1,\dots,X_N\right) := \det\left( I + \sum_{i=1}^N \alpha_i X_i \right) := \sum_{m_i} \alpha_1^{m_1}\dots\alpha_N^{m_N} Z_{m_1,\dots,m_N}\left(X_1,\dots,X_N\right) .\end{equation}

Since $\det$ is an algebraic combination of matrix entries, the expansion is a polynomial in the $\alpha_i$ variables. We can however take the sum to be over
all integer values of $m_i$ by suitably defining its $\alpha$-coefficients, $Z_{m_1,\dots,m_N}$, as zero when not corresponding to any power that appears in the $\det$ expansion.
In particular,

	\begin{equation}
	Z_{m_1,\dots,m_N}\left(X_1,\dots,X_n\right) = 0 \text{ if } \min\left\{m_i\right\}<0.
	\label{noNegativeMultiTraces}
	\end{equation}

These $Z$ coefficients then give us the \textit{generalized matrix invariants}, which reduce to the usual ones when we only consider one matrix (or set the other indices to $0$).
We can get explicit expressions in terms of traces via
\begin{equation}
\det(I+\alpha X) = e^{\mathrm{Tr}\log (I+\alpha X)}
\label{detExp}
\end{equation}
by expanding the Taylor series around $\alpha=0$ and using the linearity of the trace. This already gives, looking at the leading order of the exponential expansion,

	\begin{equation}
	Z_{0,\dots,0} (X_1,\dots,X_n) = 1.
	\label{zeroMultiTrace}
	\end{equation}

In the same way, we can reduce any expression with a $0$ index,
\begin{equation*}
Z_{n_1,\dots,n_{N-1},0} \left(X_1,\dots,X_N\right) = Z_{n_1,\dots,n_{N-1}} \left(X_1,\dots,X_{N-1}\right).
\end{equation*}

Looking at higher coefficients upon expanding the exponential returns higher invariants. For instance,
\begin{align*}
\mathrm{Tr}\log\left(I+\alpha X\right)  =& \alpha \mathrm{Tr}\left(X\right) -\frac{\alpha^2}{2}  \mathrm{Tr}\left(X^2\right) + \frac{\alpha^3}{3} \mathrm{Tr}\left(X^3\right) 
+ O\left(\alpha^4\right), \\
Z_{1}(X) =& \mathrm{Tr}(X),\\
Z_{2}(X) =& \frac{1}{2} \left( \mathrm{Tr}(X)^2-\mathrm{Tr}(X^2) \right),\\
Z_{3}(X) =& \frac{1}{6} \left( \mathrm{Tr}(X)^3-3\mathrm{Tr}(X^2)\mathrm{Tr}(X)+\mathrm{Tr}(X^3) \right),
\end{align*}
etc, but also
\begin{equation*}
Z_{1,1}(X,Y) = \mathrm{Tr}(X) \mathrm{Tr}(Y)-\mathrm{Tr}(X Y).
\end{equation*}

Of course, equation (\ref{detExp}) is usually proven using Liouville's formula. We will make contact with it again later, when looking at the derivatives of the $Z$ invariants themselves.
\subsection{Factorization}
Consider
\begin{equation*}
\det\left(1+\alpha A + \sum_{i} \beta_i B\right),
\end{equation*}
and the fact that
\begin{equation*}
\det\left( \alpha A \right) =\alpha^n \det A.
\end{equation*}
Extract this determinant from the original expansion,
\begin{align*}
&\det\left(I+\alpha A + \sum_{i} \beta_i B\right) = \sum_{l,m_i} \alpha^{l} \left(\prod \beta^{m_i}\right) Z_{l, m_1,\dots,m_N}\left(A,B_1,\dots,B_N\right) \\ =  &\det\left(A\right) \sum_{l,m_i} \alpha^{l} \left(\prod \beta_i^{m_i}\right) Z_{n-l-\sum m_i,m_1,\dots,m_N} \left(A^{-1},A^{-1}B_1,\dots,A^{-1}B_N\right).
\end{align*}
To equate the two polynomials, we equate every coefficient and get a \textbf{duality} relationship

	\begin{equation}\label{multiFactorization}
	Z_{l, m_1,\dots,m_N}\left(A,B_1,\dots,B_N\right) = \det\left(A\right) Z_{n-l-\sum m_i,m_1,\dots,m_N} \left(A^{-1},A^{-1}B_1,\dots,A^{-1}B_N\right).\end{equation}

Already an interesting property comes from the fact that any $Z$ with negative indices must be $0$, by equation (\ref{noNegativeMultiTraces}). The dual of this statement is then
\begin{equation*}
Z_{m_1,\dots,m_N}\left(X_1,\dots,X_N\right) = 0 \text{ if } \sum_{i}m_i > n.
\end{equation*}
We will call the sum of all indices $\sum m_i$ the \textbf{order} of the trace $Z_{m_1,\dots,m_N}$. 
That an invariant of order higher than the size of the matrix is zero
reduces, as expected, to the usual property of matrix invariants when we have a single matrix, and together with expression ($\ref{noNegativeMultiTraces}$) ensures that only
a finite number of $Z$ invariants for any given set of $X_i$ is non-zero.

We can also take a dual of equation $(\ref{zeroMultiTrace})$, which is the well known
\begin{equation*}
Z_{n}\left(X\right) = \det\left(X\right)
\end{equation*}
or
\begin{equation*}
1 = \det\left(X\right) Z_{n}\left(X^{-1}\right).
\end{equation*}
Now, this statement gets interesting when we introduce more matrices. Consider the two matrix case,
\begin{equation*}
Z_{l,m} \left(A,B\right) = \det\left(A\right) Z_{n-l-m,m}\left(A^{-1},A^{-1}B\right)
\end{equation*}
and set $A=X$, $B=X Y$, and $l+m=n$,
\begin{equation*}
Z_{n-m,m}\left(X, X Y\right) = \det\left(X\right) Z_{0,m}\left(X^{-1},Y\right) = \det\left(X\right) Z_m\left(Y\right).
\end{equation*}
This, which is the generalization of
\begin{equation*}
\det\left(X Y\right) = \det\left(X\right) \det\left(Y\right),
\end{equation*}
allows us to decompose order $n$ invariants of a product into products of invariants. In particular, we get the $n=2$ expression with which
we built the ODE system.

More generally,
\begin{equation}
Z_{n-\sum m_i,m_1,\dots,m_N} \left(X,X Y_1,\dots,X Y_N\right) = \det\left(X\right) Z_{m_1,\dots,m_N}\left(Y_1,\dots,Y_N\right).
\label{determinantFactorization}
\end{equation}

%
%
%
%
%
%

\subsection{Small-$\mathbf\epsilon$ expansion}
By using expression \eqref{z} we can easily derive distributivity properties, which can be applied to calculate
\begin{align*}
& Z_{l,m_1,\dots,m_N}\left(A_1+\epsilon A_2+O\left(\epsilon^2\right),B_1,\dots,B_N\right)   \\
=&\sum_{i=0}^{l} \epsilon^i Z_{l-i,i,m_1,\dots,m_N}\left(A_1,A_2+O\left(\epsilon\right),B_1,\dots,B_N\right)  \\
=&\sum_{i=0}^{l} \epsilon^i \left[Z_{l-i,i,m_1,\dots,m_N}\left(A_1,A_2,B_1,\dots,B_N\right) +O\left(\epsilon\right) \right] \\
=&\sum_{i=0}^{1} \epsilon^i Z_{l-i,i,m_1,\dots,m_N}\left(A_1,A_2,B_1,\dots,B_N\right) +O\left(\epsilon^2\right)  \\
=&Z_{l,m_1,\dots,m_N}\left(A_1,B_1,\dots,B_N\right) + \epsilon   Z_{l-1,1,m_1,\dots,m_N}\left(A_1,A_2,B_1,\dots,B_N\right) +O\left(\epsilon^2\right).
\end{align*}

\subsection{Derivatives}	
As we have seen before, the derivatives of the invariants play an essential role in the theory. We would now like to have a formula for derivatives of the form
\begin{equation*}
\frac{\mathrm{d}}{\mathrm{d}t} Z_m\left(X\left(t\right)\right).
\end{equation*}
Consider
\begin{equation*}
Z^{(m_0,m_1,m_2,\dots)}\left(X\right) := Z_{m_0,m_1,m_2,\dots}\left(X,X',X'',\dots\right),
\end{equation*}
such that for some $N$ we have $m_i=0$ for every $i>N$.

We can retrieve its first derivative from its Taylor series, which we can in turn get from its small $\epsilon$ expansion.
\begin{align*}
& Z^{(m_0,m_1,\dots)}\left( X+\epsilon X'+O\left(\epsilon^2\right) \right)  \\
=&Z^{(m_0,m_1,m_2,\dots)}\left(X\right)+\epsilon   \left(m_1+1\right) Z_{m_0-1,m_1+1,m_2,\dots}\left(X,X',X'',\dots\right) \\ & +\epsilon   \left(m_2+1\right) Z_{m_0,m_1-1,m_2+1,\dots}\left(X,X',X'',\dots\right)+\cdots
\end{align*}

Taking the $\epsilon$ term we get the first coefficient of the Taylor series, i.\,e., the first derivative,

\[
	\left(Z^{(m_0,m_1,m_2,\dots)}\left(X\right)\right)' =\sum_{i=1}^\infty \left(m_i+1\right) Z^{ (m_0,\dots,m_{i-1}-1,m_i+1,\dots    ) }\left(X\right).
	\]
Notice that the infinite sum is merely formal, since by equation (\ref{noNegativeMultiTraces}) it is guaranteed to terminate as soon as all the remaining $m_i$ are $0$, due
to the $m_{i-1}-1$ index at every term.


For the first few derivatives, we find via recursion the general expressions
\begin{align*}
Z_{m} \left(X\right)' = \left(Z^{(m)}\right)'=&Z^{(m-1,1)},
\\
Z_{m} \left(X\right)''= \left(Z^{(m)}\right)''=& 2 Z^{(m-2,2)} + Z^{(m-1,0,1)},
\\
Z_{m} \left(X\right)''' = \left(Z^{(m)}\right)''' =& 6 Z^{(m-3,3)} + 3 Z^{(m-2,1,1)} + Z^{(m-1,0,0,1)}.
\end{align*}

Something very important (albeit somehow obvious, following Leibniz's rule for matrices), is that the order of the invariants involved in the expressions is preserved.


This allows us to use the factorization formula (\ref{determinantFactorization}) over the derivatives of the determinant, which corresponds to $Z_n$.

As a small note, if we take in the first derivative $m=n$ together with expression (\ref{multiFactorization}), we get
\begin{equation*}
\det\left(X\right)'=Z_{n-1,1}\left(X,X'\right) = \det\left(X\right) Z_1\left(X^{-1} X'\right) = \det\left(X\right) \mathrm{Tr}\left(X^{-1} X'\right),
\end{equation*}
the usual Liouville's Formula.

\subsection{Application to the differential system of invariants for $\mathbf{n>2}$}

In the matrix dimension $m=2$ case, taking derivatives of the determinant eventually closes, since $X''=XE$. This follows from
\[\det(X)'' = Z_m(X)'' = 2 Z^{(m-2,2)}(X) + Z^{(m-1,0,1)}(X).\]
The $Z^{(m-1,0,1)}(X)$ can be immediately rewritten as a determinant by using the duality formula,
\[Z_{m-1,0,1} (X, X', X'') = \det(X) Z_1 (X^{-1} X'') = \det(X) \tr(E),\]
and, when $m=2$,
\[Z^{(m-2,2)}(X) = Z^{(0,2)}(X) = \det(X').\]
Of course, now we can do the same for $X'$,
\[\det(X')'' = 2 Z^{(m-2,2)}(X) + Z^{(m-1,0,1)}(X)\]
and, for $m=2$,
\[\det(X')'' = 2 \det(X'') + \det(X') \tr(E) = 2 \det(E) \det(X) + \det(X') \tr(E),\]
closing the system as we had found in Theorem \ref{AJLid}. The problem is now obvious, since for $m>2$, $Z^{(m-2,2)}(X)$ will involve a non trivial product between $X$ and $X'$. One could consider this as a new variable for the system, but its derivatives will now concern objects of the form $Z^{(m-2,1,1)}(X)$ which, if understood as yet another variable of the system, would yield upon derivation
\[Z^{(m-2,1,0,1)}(X),\ Z^{(m-2,1,0,0,1)}(X),\ Z^{(m-2,1,0,...,0,1)}, ...\]

Notice that this will always involve a term in $X'$, and a term in $X$, so that we cannot perform the same trick as we did for $Z^{(m-1,0,1)}(X)$ --namely, using $X''=XE$ to factor the determinant out. Hence, the system of second derivatives of invariants for $m>2$ does not close.

\end{document}